\magnification=\magstep1

          %%%%%%%%%%%%%%%%%%%%%%%%%%%%%%%%%%%%%%%%%%%%%%%%%%%%%%%%%%%
          %       FOP.tex FORTEX-compatible version of papmacs      %
          %   VERSION OF May 6, 1989; does not use the AMS-fonts    %
          %%%%%%%%%%%%%%%%%%%%%%%%%%%%%%%%%%%%%%%%%%%%%%%%%%%%%%%%%%%

\def\item{\vskip1.3pt\hang\textindent}% THIS REPLACES KNUTH'S DEF'N

% THIS
                                              %REPLACES KNUTH'S DEF'N

\tolerance=300 \pretolerance=200 \hfuzz=1pt \vfuzz=1pt

%\magnification=\magstep1
\hoffset 0cm            %all offsets went wrong
%\voffset=0.8true cm
\hsize=5.8 true in \vsize=9.5 true in

\def\rightheadline{\hfil\smc\lastname\hfil\tenbf\folio}
\def\leftheadline{\tenbf\folio\hfil\smc\lastname\hfil}
\headline={\ifodd\pageno\rightheadline\else\leftheadline\fi}
\newdimen\dimenone
\def\checkleftspace#1#2#3#4#5{%DIESER MACRO STAMMT VON APPELT
 \dimenone=\pagetotal
 \advance\dimenone by -\pageshrink   %testen ob Titel noch mit Gewalt auf Seite
                                                                          %geht
 \ifdim\dimenone>\pagegoal          %nacha tua nix-- gewoehnliche Outputroutine
   \else\dimenone=\pagetotal
        \advance\dimenone by \pagestretch
        \ifdim\dimenone<\pagegoal
          \dimenone=\pagetotal
          \advance\dimenone by#1         %addieren Skip vor Ueberschrift (=#1)
          \setbox0=\vbox{#2\parskip=0pt                %#2 ist gewaehlter Font
                       \hyphenpenalty=10000
                       \rightskip=0pt plus 5em
                       \noindent#3 \vskip#4}    %#3=Ueberschrift,#4=skip nachher
        \advance\dimenone by\ht0
        \advance\dimenone by 3\baselineskip
        \ifdim\dimenone>\pagegoal\vfill\eject\fi
          \else\eject\fi\fi}

\parindent=35pt
\mathsurround=1pt
\parskip=1pt plus .25pt minus .25pt
\normallineskiplimit=.99pt

\mathchardef\emptyset="001F % THIS REPLACES KNUTH'S DEFINITION

% USED FOR REAL PART OF COMPLEX NUMBERS
% USED FOR IMAGINARY PART OF COMPLEX NUMBERS
 % USED FOR IDENTITY FUNCTION
\def\Int{\mathop{\rm int}\nolimits}
% USED FOR TRACE OF MATRIX
%

%\def\u{\mathop{\rm u}\nolimits}

%\def\O{\mathop{\rm O}\nolimits}

\def\1{{\bf1}}\def\0{{\bf0}}

\def\({\bigl(}  \def\){\bigr)}
\def\<{\mathopen{\langle}}\def\>{\mathclose{\rangle}}

\def\Z{{\mathchoice{{\hbox{$\rm Z\hskip 0.26em\llap{\rm Z}$}}}%
{{\hbox{$\rm Z\hskip 0.26em\llap{\rm Z}$}}}%
{{\hbox{$\scriptstyle\rm Z\hskip 0.31em\llap{$\scriptstyle\rm Z$}$}}}{{%
\hbox{$\scriptscriptstyle\rm
Z$\hskip0.18em\llap{$\scriptscriptstyle\rm Z$}}}}}}

\def\F{{\mathchoice{\hbox{$\rm I\hskip-0.14em F$}}%
{\hbox{$\rm I\hskip-0.14em F$}}%
{\hbox{$\scriptstyle\rm I\hskip-0.14em F$}}%
{\hbox{$\scriptscriptstyle\rm I\hskip-0.10em F$}}}}

\def\R{{\mathchoice{\hbox{$\rm I\hskip-0.14em R$}}%
{\hbox{$\rm I\hskip-0.14em R$}}%
{\hbox{$\scriptstyle\rm I\hskip-0.14em R$}}%
{\hbox{$\scriptscriptstyle\rm I\hskip-0.10em R$}}}}

\def\K{{\mathchoice{\hbox{$\rm I\hskip-0.15em K$}}%
{\hbox{$\rm I\hskip-0.15em K$}}%
{\hbox{$\scriptstyle\rm I\hskip-0.15em K$}}%
{\hbox{$\scriptscriptstyle\rm I\hskip-0.11em K$}}}}

\def\qed{\hfill {\hbox{[\hskip-0.05em ]}}}

\def\.{{\cdot}}
\def\|{\Vert}
\def\ssk{\smallskip}
\def\msk{\medskip}
\def\bsk{\bigskip}
\def\giantskip{\vskip2\bigskipamount}

\def\giantbreak{\par \ifdim\lastskip<2\bigskipamount \removelastskip
         \penalty-400 \giantskip\fi}

\def\nin{\noindent}
\def\cen{\centerline}
\def\pagebreak{\vskip 0pt plus 0.0001fil\break}
\def\linebreak{\break}

\def\epsilon{\varepsilon}

\font\ninerm=cmr9 \font\eightrm=cmr8 \font\sixrm=cmr6
 \font\eightbf=cmbx8 \font\sixbf=cmbx6
 \font\eighti=cmmi8 \font\sixi=cmmi6
\font\ninesy=cmsy9 \font\eightsy=cmsy8 \font\sixsy=cmsy6
 \font\eightit=cmti8 
% SANS SERIF 10 POINT
 %SANS SERIF 10 POINT ITALIC
 \font\eightsl=cmsl8 
\font\eighttt=cmtt8
 %SLANTED TYPEWRITER 10 POINT
 %BOLD FACE MATH SYMBOLS 10 POINT
 %DUNHILL STYLE 10 POINT
 %SAN SERIF BOLD EXTENDED 10 POINT
 %USED FOR TITLES
 %USED FOR TITLES
\font\bfone=cmbx10 scaled\magstep1 %BOLDFACE AT MAGSTEP 1
 %BOLDFACE AT MAGSTEP 2
 %BOLDFACE AT MAGSTEP 3
\font\smc=cmcsc10 
 
scaled\magstep1 \font\small=cmcsc8

\def\no #1. {\bigbreak\vskip-\parskip\noindent\bf #1. \quad\rm}

\def\Proposition #1. {\checkleftspace{0pt}{\bf}{Theorem}{0pt}{}
\bigbreak\vskip-\parskip\noindent{\bf Proposition #1.} \quad\it}

\def\Theorem #1. {\checkleftspace{0pt}{\bf}{Theorem}{0pt}{}
\bigbreak\vskip-\parskip\noindent{\bf  Theorem #1.} \quad\it}
\def\Corollary #1. {\checkleftspace{0pt}{\bf}{Theorem}{0pt}{}
\bigbreak\vskip-\parskip\nin{\bf Corollary #1.} \quad\it}
\def\Lemma #1. {\checkleftspace{0pt}{\bf}{Theorem}{0pt}{}
\bigbreak\vskip-\parskip\noindent{\bf  Lemma #1.}\quad\it}

\def\Definition #1. {\checkleftspace{0pt}{\bf}{Theorem}{0pt}{}
\rm\bigbreak\vskip-\parskip\noindent{\bf Definition #1.} \quad}

\def\Remark #1. {\checkleftspace{0pt}{\bf}{Theorem}{0pt}{}
\rm\bigbreak\vskip-\parskip\noindent{\bf Remark #1.}\quad}

\def\Exercise #1. {\checkleftspace{0pt}{\bf}{Theorem}{0pt}{}
\rm\bigbreak\vskip-\parskip\noindent{\bf Exercise #1.} \quad}

\def\Example #1. {\checkleftspace{0pt}{\bf}{Theorem}{0pt}{}
\rm\bigbreak\vskip-\parskip\noindent{\bf Example #1.}\quad}
\def\Examples #1. {\checkleftspace{0pt}{\bf}{Theorem}{0pt}
\rm\bigbreak\vskip-\parskip\noindent{\bf Examples #1.}\quad}

\newcount\problemnumb \problemnumb=0
\def\Problem{\global\advance\problemnumb by 1\bigbreak\vskip-\parskip\noindent
{\bf Problem \the\problemnumb.}\quad\rm }

\def\Proof#1.{\rm\par\ifdim\lastskip<\bigskipamount\removelastskip\fi\smallskip
            \noindent {\bf Proof.}\quad}

\nopagenumbers

\def\author{}
\def\lastname{}
\def\thanks#1{\footnote*{\eightrm#1}}
\def\title{}

\def\lastname{}
\def\h{{\textstyle{1\over2}}}

\def\n{{\cal N}}
\def\ep{\epsilon}

\def\text{\textstyle}
\def\disp{\displaystyle}
\def\d{{\,\rm d}}

\def\and{{\rm and }}

\def\n{\cen{{\it W.G. Nowak}}}

\expandafter\edef\csname amssym.def\endcsname{%
       \catcode`\noexpand\@=\the\catcode`\@\space}
%  Set the catcode to 11 for use in private control sequence names.
\catcode`\@=11
\def\undefine#1{\let#1\undefined}
\def\newsymbol#1#2#3#4#5{\let\next@\relax
 \ifnum#2=\@ne\let\next@\msafam@\else
 \ifnum#2=\tw@\let\next@\msbfam@\fi\fi
 \mathchardef#1="#3\next@#4#5}
\def\mathhexbox@#1#2#3{\relax
 \ifmmode\mathpalette{}{\m@th\mathchar"#1#2#3}%
 \else\leavevmode\hbox{$\m@th\mathchar"#1#2#3$}\fi}
\def\hexnumber@#1{\ifcase#1 0\or 1\or 2\or 3\or 4\or 5\or 6\or 7\or 8\or
 9\or A\or B\or C\or D\or E\or F\fi}

%Loading fontfiles `eufm' and `msbm'

\font\tenmsb=msbm10 \font\sevenmsb=msbm7 \font\fivemsb=msbm5
\newfam\msbfam
\textfont\msbfam=\tenmsb\scriptfont\msbfam=\sevenmsb
\scriptscriptfont\msbfam=\fivemsb \edef\msbfam@{\hexnumber@\msbfam}
\def\Bbb#1{{\fam\msbfam\relax#1}}

\font\teneufm=eufm10 \font\seveneufm=eufm7 \font\fiveeufm=eufm5
\newfam\eufmfam
\textfont\eufmfam=\teneufm \scriptfont\eufmfam=\seveneufm
\scriptscriptfont\eufmfam=\fiveeufm

\catcode`@=11 %Set the catcode to 11 for use in private control sequence names.

\expandafter\edef\csname amssym.def\endcsname{%
       \catcode`\noexpand\@=\the\catcode`\@\space}
\font\eightmsb=msbm8 \font\sixmsb=msbm6 \font\fivemsb=msbm5
\font\eighteufm=eufm8 \font\sixeufm=eufm6 \font\fiveeufm=eufm5
\newskip\ttglue
\def\eightpoint{\def\rm{\fam0\eightrm}%
  \textfont0=\eightrm \scriptfont0=\sixrm \scriptscriptfont0=\fiverm
  \textfont1=\eighti \scriptfont1=\sixi \scriptscriptfont1=\fivei
  \textfont2=\eightsy \scriptfont2=\sixsy \scriptscriptfont2=\fivesy
  \textfont3=\tenex \scriptfont3=\tenex \scriptscriptfont3=\tenex
\textfont\eufmfam=\eighteufm \scriptfont\eufmfam=\sixeufm
\scriptscriptfont\eufmfam=\fiveeufm \textfont\msbfam=\eightmsb
\scriptfont\msbfam=\sixmsb \scriptscriptfont\msbfam=\fivemsb
  \def\it{\fam\itfam\eightit}%
  \textfont\itfam=\eightit
  \def\sl{\fam\slfam\eightsl}%
  \textfont\slfam=\eightsl
  \def\bf{\fam\bffam\eightbf}%
  \textfont\bffam=\eightbf \scriptfont\bffam=\sixbf
   \scriptscriptfont\bffam=\fivebf
  \def\tt{\fam\ttfam\eighttt}%
  \textfont\ttfam=\eighttt
  \tt \ttglue=.5em plus.25em minus.15em
  \normalbaselineskip=9pt
  \def\MF{{\manual opqr}\-{\manual stuq}}%
  \let\big=\eightbig
  \setbox\strutbox=\hbox{\vrule height7pt depth2pt width\z@}%
  \normalbaselines\rm}
\def\eightbig#1{{\hbox{$\textfont0=\ninerm\textfont2=\ninesy
  \left#1\vbox to6.5pt{}\right.\n@space$}}}

%  Restore the catcode value for @ that was previously saved.

\csname amssym.def\endcsname

% end of \symb1.tex

\def\la{\lambda}
\def\al{\alpha}
\def\be{\beta}

\def\({\left(}
\def\){\right)}
\def\for{\qquad \hbox{for}\ }
\def\eq{\eqalign}

\def\O#1{O\(#1\)}
\def\abs#1{\left| #1 \right|}

\def\klein{\eightpoint \def\smc{\small} \baselineskip=9pt}

\def\fn#1#2{{\parindent=0.7true cm
\footnote{$^{(#1)}$}{{\klein  #2}}}}

\font\boldmas=msbm10                  %%
\def\Bbb#1{\hbox{\boldmas #1}}        %%
\def\Z{{\Bbb Z}}                        %%

\def\R{{\Bbb R}}
\def\F{{\Bbb F}}

                  %%
        %%

%%%%%%%%%%%%%%%%%%%%%%%%%%%%%%%%%%%%%%%%%%%%%%%%%%%%%%%%%%%%%%%%%%%%%%
\font\eightrm=cmr8 \long\def\fussnote#1#2{{\baselineskip=9pt
\setbox\strutbox=\hbox{\vrule height 7pt depth 2pt width 0pt}%
\eightrm \footnote{#1}{#2}}}
%%%%%%%%%%%%%%%%%%%%%%%%%%%%%%%%%%%%%%%%
%% zum Verkleinern (Summationsgrenzen, Folgenindex, e.a.)%%
%%%%%%%%%%%%%%%%%%%%%%%%%%%%%%%%%%%%%%%%
\font\boldmasi=msbm10 scaled 700      %%
\def\Bbbi#1{\hbox{\boldmasi #1}}      %%
\font\boldmas=msbm10                  %%
\def\Bbb#1{\hbox{\boldmas #1}}        %%
\def\Pi{{\Bbbi P}}                      %%
\def\Ri{{\Bbbi R}}

                      %%
%%%%%%%%%%%%%%%%%%%%%%%%%%%%%%%%%%%%%%%%

                        %%

\def\dint #1 {% Doppelintegral; #1 Integrationsbereich
\quad  \setbox0=\hbox{$\disp\int\!\!\!\int$}
  \setbox1=\hbox{$\!\!\!_{#1}$}
  \vtop{\hsize=\wd1\centerline{\copy0}\copy1} \quad}

\def\drint #1 {% Dreifachintegral; #1 Integrationsbereich
\qquad  \setbox0=\hbox{$\disp\int\!\!\!\int\!\!\!\int$}
  \setbox1=\hbox{$\!\!\!_{#1}$}
  \vtop{\hsize=\wd1\centerline{\copy0}\copy1}\qquad}

\def\frac#1#2{{#1\over #2}}

\def\date{\the\day.~\the\month.~\the\year}

\def\klein{\eightpoint \def\smc{\small} }

\def\frac#1#2{{#1\over#2}}
\def\Int{\int\limits}

\def\vol{{\rm vol}}

%\nonumbers
\hsize=16true cm     \vsize=23true cm

\parindent=0cm

\def\J{{\cal J}}
\def\K{{\cal K}}
\def\r{{\cal R}}
\def\c{{\cal C}}
\def\B{{\cal B}}
\def\D{{\cal D}}
\def\F{{\cal F}}
\def\dB{\partial\B}
\def\dr{\partial\r}
\def\InT{\Int_T^{2T}}
\def\kn{(k/\sqrt{n})}
\def\n{\sqrt{n}}
\def\ni{\sqrt{n'}}
\def\fou#1{\widehat{#1}}
\def\b#1{{\bf #1}}

\def\lT{(L/T)}

\vbox{\vskip 1.2true cm}

\cen{{\bfone A mean-square bound for the lattice discrepancy of}}
\msk \cen{{\bfone bodies of rotation with flat points on the
boundary }}\bsk \cen{{\bf Werner Georg Nowak}\fn{*}{The author
gratefully acknowledges support from the Austrian Science Fund (FWF)
under project Nr.~P18079-N12.} {\bf(Vienna)}}

\vbox{\vskip 0.7true cm}

\footnote{}{\klein{\it Mathematics Subject Classification }
(2000): 11P21, 11K38, 52C07.\par }

\cen{\it Dedicated to Professor Edmund Hlawka on his 90th birthday}

%{\klein{\bf Abstract. } }

\vbox{\vskip 1true cm}

{\bf 1.~Introduction.}\quad Let $\B$ denote a compact convex body in
$s$-dimensional Euclidean space, $s\ge2$, which contains the origin
as an inner point and whose boundary $\dB$ is sufficiently smooth.
The central question of the classic {\it lattice point theory of
large domains } consists of estimating the {lattice point
discrepancy} of a linearly dilated copy $t\,\B$, i.e.,
$$ P_\B(t) := \#\(t\,\B\cap\Z^s\) - \vol(\B)t^s\,, \eqno(1.1) $$
where $t$ is a large real parameter. For enlightening accounts on
this topic, the reader is referred to E.~Kr\"atzel's monographs [14]
and [15], to a recent survey article by A.~Ivi\'c, E.~Kr\"{a}tzel,
M.~K\"{u}hleitner, and W.G.~Nowak [10], and to M.~Huxley's book [7]
where he exposed his breakthrough in the planar case ({"Discrete
Hardy-Littlewood method"}). \ssk For the case that $\dB$ is of
bounded nonzero Gaussian curvature throughout, the usual and
plausible conjecture is that
$$ P_\B(t) \ll t^{\Theta_s+\ep} \eqno(1.2)$$
for every $\ep>0$, where $\Theta_2:=\h$, $\Theta_s:=s-2$ for
$s\ge3$. It is well-known that for every dimension, $\Theta_s$ is
the minimal possible value with this property, and that (1.2) is
actually true for spheres of dimension $s\ge4$, even with $\ep=0$ if
$s\ge5$: See, e.g., E.~Kr\"atzel [15], p.~227. Quite recently,
V.~Bentkus and F.~G\"otze [1] and F.~G\"otze [3] established (1.2)
for arbitrary ellipsoids of dimension $s\ge5$. \ssk However, for
$s=2$ and $3$, and for general bodies of higher dimensions, the
proof or disproof of (1.2) remains an open problem. The sharpest
known estimates are due to W.~M\"{u}ller [24]. Researchers subsequently
dealt with the task to verify (1.2) "on average", i.e., to show that
$$ \Int_0^T \(P_\B(t)\)^2 \d t \ll T^{2\Theta_s+1+\ep}\,. \eqno(1.3)
$$ In fact, (1.3) was established for planar domains by D.G.~Kendall [13] and the
author [25], ultimately in the form of an asymptotics [26]. For
dimensions $s\ge4$, (1.3) was proved by W.~M\"uller [23] who however
had to leave open the case $s=3$. This gap was filled by
A.~Iosevich, E.~Sawyer, \and A.~Seeger [8] who showed that $$
\Int_0^T \(P_\B(t)\)^2 \d t \ll \cases{T^{2\Theta_s+1} & for
$s\ge4$,\cr T^3\,(\log T)^2 & for $s=3$.\cr }\eqno(1.4)$$ The very
last estimate comes rather close to the asymptotic formula known for
the three-dimensional sphere $\B_0$, namely $$\Int_0^T
\(P_{\B_0}(t)\)^2 \d t = C\,T\,\log T + \O{T\,(\log T)^{1/2}}\,. $$
Cf.~V.~Jarnik [11], and also Y.-K.~Lau [22] who improved the error
term to $O(T)$. \bsk\msk

{\bf 2.~Recent developments and statement of present result. } The
topic of this note will combine two recent trends in lattice point
theory: On the one hand, increased interest arose in $\R^3$-bodies
{\it of rotation} (with respect to one of the coordinate axes),
denoted by $\r$ in what follows. For the case of nonzero curvature,
F.~Chamizo [2] obtained the upper bound
$$ P_\r(t)\ll t^{11/8+\ep}\,, \eqno(2.1)
$$ while papers by M.~K\"{u}hleitner [20] and M.~K\"{u}hleitner and
W.G.~Nowak [21] provided $\Omega$-results. A recent article of
E.~Kr\"{a}tzel and W.G.~Nowak [19] gives a version of (2.1) with
numerical constants, for the special case of an ellipsoid. \ssk

On the other hand, a number of papers investigated the influence of
boundary points with curvature zero on the lattice discrepancy.
While Kr\"{a}tzel's monograph [15] provides an enlightening survey on
the planar case (which is comparatively well understood), results
for dimension 3 and higher can be found in the works of K.~Haberland
[4], E.~Kr\"{a}tzel [16], [17], [18], and M.~Peter [27]. These are all
(pointwise) $O$-estimates, partially providing a precise evaluation
of the contribution of an isolated flat point on $\dr$ to $P_\r(t)$,
with a remainder of smaller order. \msk

In the present paper we shall take up both of these matters, under
the aspect of a mean-square estimate in the sense of (1.3). We
will consider a $\R^3$-body of rotation $\r$ (with respect to one
of the coordinate axes), with smooth boundary $\partial\r$ of
nonzero Gaussian curvature $\kappa$ throughout, except for the
points of intersection of $\partial\r$ with the axis of rotation,
where $\kappa$ may vanish. It will turn out that the contribution
of these flat points to the lattice point discrepancy can be
evaluated quite accurately, leaving a remainder term which is in
mean-square "as small as it should be", in the sense of formula
(1.3). \ssk We remark parenthetically that if $\kappa$ would
vanish anywhere else on $\dB$, it would do so on a whole circle.
This will presumably have a more dramatic effect on the lattice
discrepancy. It seems much more difficult to obtain a sharp result
in this general case. \ssk

{\it Precise formulation of present assumptions. } Let $\rho:
[0,\pi] \rightarrow \R_{>0}$ be a function of class $C^4$, with
$\rho'(0)=\rho'(\pi)=0$ and\fn{1}{Recall that the curvature of a
curve whose equation in polar coordinates is $\rho=\rho(\theta)$,
is given, in absolute value, by
${\abs{\rho(\theta)\,\rho''(\theta)-2
\rho'^2(\theta)-\rho^2(\theta)}
/(\rho^2(\theta)+\rho'^2(\theta))^{3/2}}$.}
$$ \rho\,\rho''-2 \rho'^2-\rho^2\ne0 \eqno(2.2)$$ throughout
$]0,\pi[$. Suppose that $\rho$ is analytic in $\pi$ and 0. At
these two values, the left-hand side of (2.2) may vanish, of
orders (exactly) $N_1, N_2\ge0$, as a function of $\theta$, the
case that $\min(N_1,N_2)=0$ not being excluded. Then
$$ \c = \{ (x,y)=(\rho(|\theta|)\cos\theta,
\rho(|\theta|)\sin\theta):\ \theta\in[-\pi,\pi]\ \} $$ defines a
smooth curve in the $(x,y)$-plane, symmetric with respect to the
$x$-axis. Rotating $\c$ around the latter, we obtain a smooth
surface in $(x,y,z)$-space, which we call $\partial\B$, where $\B$
is the compact convex body bounded by $\partial\B$. We denote by
$a_1, a_2$ the minimal, resp., maximal $x$-coordinate on
$\partial\B$. Obviously, the Gaussian curvature of $\partial\B$
vanishes at most in the points of intersection with the $x$-axis.
\bsk

{\bf Theorem. } {\it Suppose that the conditions stated are
satisfied, in particular, that the Gaussian curvature of $\dr$
vanishes at most in the two points of intersection with the axis
of rotation. Then for the number $A_\r(t)$ of lattice points in
the linearly dilated body $t\,\r$ the asymptotic formula
$$ \eq{A_\r(t) &= \vol(\r)\,t^3 + \sum_{j=2}^{N_1+1} d_{1,j}^*\,
\F(-a_1t,j/(N_1+2))\,t^{2-j/(N_1+2)} - \cr &- \sum_{j=2}^{N_2+1}
d_{2,j}^*\, \F(a_2t,j/(N_2+2))\,t^{2-j/(N_2+2)} +
\Delta_\r(t)\,,\cr } $$ holds true, where $$
\F(\xi,\eta):=(2\pi)^\eta\,\Gamma(\eta)\,\sum_{k=1}^\infty
k^{-1-\eta}\,\sin(2\pi k\xi - \h\pi\eta) \qquad(\eta>0)\,, $$ and
the remainder satisfies the mean-square estimate $$\Int_0^T
\(\Delta_\r(t)\)^2\,\d t = \O{T^{3+\ep}}$$ for each $\ep>0$. The
coefficients $d_{1,j}^*$, $d_{2,j}^*$ are computable, on the basis
of the formulas $(3.1)$ - $(3.5)$ below. In particular,
$d_{1,1}^*>0$, $d_{2,1}^* <0$.} \bsk

{\bf Remarks. } 1. It is easy to see that the error term satisfies
in fact the pointwise upper estimate $$
\Delta_\r(t)=\O{t^{3/2+\ep}}\,. \eqno(2.3)$$ This is a
straightforward consequence of the works of Kr\"{a}tzel [15], [16],
[17], [18], but follows also as a simple by-result from the argument
in this paper: See the concluding remark at the end. \msk

2. There is a crucial difference in the treatment of the problem,
depending on whether there are boundary points of curvature zero or
not. For $\kappa$ nonzero, the analysis leading to the results
(1.3), (1.4) is based on the asymptotic expansion of the Fourier
transform of the indicator function of the body $\B$, which is due
to E.~Hlawka [5], [6]. In the case that $\inf\kappa=0$ the latter is
not at our disposal. Thus we have to employ a quite different
approach which uses a truncated Hardy's identity (Lemma 1) and a
transformation of trigonometric sums.

\bsk\msk

{\bf 3.~Some auxiliary results. } \bsk

{\bf Lemma 1. } { \it For integers $k\ge0$, let as usual $r(k)$
denote the number of pairs $(m_1, m_2)\in\Z^2$ with $m_1^2+m_2^2=k$.
For large real parameters $X, Y$ with $ \log Y\ll \log X$, and any
$\ep>0$, it then follows that $$ \eq{P(X)&:=\sum_{0\le k\le
X}r(k)-\pi X = \cr & = {1\over\pi}\,X^{1/4}\sum_{1\le n\le
Y}{r(n)\over n^{3/4}}\,\cos(2\pi\sqrt{nX}-3\pi/4)\ +
\O{X^{1/2+\ep}\,Y^{-1/2}}+\O{Y^\ep}\,.\cr} $$} \bsk

{\bf Proof. } This is contained in formula (1.9) of A.~Ivi\'c [9].
\qed \bsk\msk

{\bf Lemma 2A. } { \it Let $F\in C^4[A, B]$,  $G\in C^2[A, B]$,
and suppose that, for positive parameters $X, Y, Z$, we have
$B-A\ll X$ and
$$ F^{(j)}\ll X^{2-j} Y^{-1} \for j=2,3,4, \ \
F''\ge c_0 Y^{-1}\,,\quad G^{(j)}\ll X^{-j} Z \for j=0,1,2, $$
throughout the interval $[A,B]$, with some constant $c_0>0$.
Assume further that there exists a value $\tau_0\in]A,B[$ with
$F'(\tau_0)=0$. Writing as usual $e(z):=e^{2\pi{\rm i}z}$, it
follows that
$$ \eq{\Int_A^B G(\tau) e(F(\tau))\d\tau =&\
{G(\tau_0)\over\sqrt{F''(\tau_0)}}\,e(F(\tau_0)+1/8) + \O{X^{-1}YZ}
+ \cr & + \O{Z\,\min\(|F'(A)|^{-1},\sqrt{Y}\)} +
\O{Z\,\min\(|F'(B)|^{-1},\sqrt{Y}\)}\,. \cr } $$ } \bsk

{\bf Proof. } This is Lemma III.2 in A.A.~Karatsuba and S.M.~Voronin
[12].    \qed \bsk

{\bf Lemma 2B. } { \it Let $F\in C^4[A, B]$,  $G\in C^2[A, B]$, and
suppose that, for positive parameters $X, Y, Z$, we have $1\le
B-A\ll X$ and
$$ F^{(j)}\ll X^{2-j} Y^{-1} \for j=2,3,4, \
\abs{F''}\ge c_0 Y^{-1}\,,\quad G^{(j)}\ll X^{-j} Z \for j=0,1,2,
$$ throughout the interval $[A,B]$, with some constant $c_0>0$.
Let $\J'$ denote the image of $[A,B]$ under $F'$, and $F^*$ the
inverse function of $F'$. Then
$$ \eq{\sum_{A<m\le B} G(m)\,e(F(m)) =&\ e\({{\rm sgn}(F'')\over8}\) \sum_{k\in\J'}{G(F^*(k))
\over\sqrt{\abs{F''(F^*(k))}}}\,e\(F(F^*(k))-kF^*(k)\)  + \cr & +
\O{Z\(\sqrt{Y}+1+{Y\over X}+\sum_{\Lambda=X,Y,Z}|\log\Lambda|\)}\,.
\cr }
$$ }\bsk

{\bf Proof. } This is essentially contained in Theorem 2.11 of
E.~Kr\"{a}tzel's monograph [14], apart from his cumbersome condition
(2.37) which basically requires the function $F$ to be algebraic.
But this can be avoided by replacing, in Kr\"{a}tzel's proof, his Lemma
2.5 by the result we just stated as Lemma 2A. \qed

\bsk

For our argument it will be essential to have at hand a close
analysis of the situation near the points were the Gaussian
curvature (possibly) vanishes. To this end, let $\c^+$ denote the
upper half of $\c$, and set
$$ \c^+= \{(\rho(\theta)\cos\theta,
\rho(\theta)\sin\theta):\ \theta\in[0,\pi]\ \} = \{(x,y):\ a_1\le
x\le a_2\,,\ y=f(x)\ \}\,,   $$ with
$a_1:=-\rho(\pi)<0<a_2:=\rho(0)$. This defines $f:\ [a_1,
a_2]\rightarrow\R_{\ge0}$ as a strictly positive $C^4$-function on
$]a_1,a_2[$, with $f(a_1)=f(a_2)=0$, and $f''$ strictly negative
throughout. By our assumptions, for each of the $a_i$'s, and
$(x,y)\in\c^+$ in a suitable neighborhood of $(a_i,0)$, $$ x = a_i +
c_i\,y^{N_i+2}+\sum_{m=1}^\infty
c_{i,m}\,y^{N_i+2+m}\qquad(c_i\ne0)\,.\eqno(3.1)
$$ Consequently, $$ \eq{y = f(x) &= \sum_{j=1}^\infty
d_{i,j}\,|x-a_i|^{\al_i\,j}\,,\cr
 \al_i &:= {1\over N_i+2}\,,\ d_{i,1}=|c_i|^{-\al_i}\ne0\,,\cr} \eqno(3.2) $$ and
the other $d_{i,j}$'s can be computed recursively from the
$c_{i,m}$'s. It thus follows that, for $r=0,1,2,\dots,$
$$ f^{(r)}(x) \asymp |x-a_i|^{\al_i-r} \eqno(3.3)  $$ for $x$ close
to $a_i$. Similarly, we deduce that $$ {\d^r\over\d
x^r}\(\sqrt{f(x)}\) \asymp |x-a_i|^{\al_i/2-r} \eqno(3.4)  $$ for
$r=0,1,2,\dots,$ and $x$ near $a_i$. Furthermore,

\vbox{$$ {\d\over\d x}\(f^2(x)\) = \sum_{j=2}^\infty d_{i,j}^*
|x-a_i|^{\al_i\,j-1}\,, \eqno(3.5) $$ with $d_{i,2}^*=2\al_i
d_{i,1}^2(-1)^{i+1}$, again in an appropriate neighborhood of
$a_i$.} \ssk For our proof we will also need some knowledge about
the {\it tac-function} of $\r$ $$
H(u,v,w):=\max_{(x,y,z)\in\r}(ux+vy+wz)
$$ and the {\it polar body } $\r^*:\ H(u,v,w)\le1$. The
connection between the respective smoothness and the curvature of
$\partial\r$ and of $\partial\r^*$ has been neatly worked out in
W.~M\"{u}ller [23], Lemma 1. It is clear that $H$ and thus $\r^*$ is
again invariant under rotations around the first coordinate axis.
Let $C_+^*$ denote the intersection of $\partial\r^*$ with the
closed upper half of the $(u,v)$-plane. Then $$ C_+^* = \{\ (u,v):\
{1/a_1}\le u \le 1/a_2\,,\ v = h(u) \}\,,  $$ where $h:\ [1/a_1,
1/a_2]\rightarrow\R_{\ge0}$ is a strictly positive $C^3$-function on
$]1/a_1,1/a_2[$, with $h(1/a_1)=h(1/a_2)=0$ (cf.~W.~M\"{u}ller [23],
Lemma 1). \bsk

{\bf Lemma 3A. } { \it With the conditions and definitions stated,
$$ \sup_{1/a_1<u<1/a_2}\,|h(u)\,h'(u)| < \infty\,. $$ } \bsk

{\bf Proof. } Let the real numbers $x\in\ ]a_1,a_2[$ and $u\in\
]1/a_1,1/a_2[$ be connected by the condition\fn{2}{In fact, the
points $(x,f(x))$ and $(u,h(u))$ are called {\it polar reciprocal}
to each other. The correspondence defined is one-one.}
$$ 1=H(u,h(u),0)=\max_{\xi\in[a_1,a_2]}
(u\xi+h(u)f(\xi))=ux+h(u)f(x)\,. $$ Plainly, $x\to a_i$ if
$u\to1/a_i$, and vice versa. Eliminating $h(u)$ from the pair of
equations $$ \eq{ux+h(u)f(x)&=1\,,\cr u+h(u)f'(x)&=0\,,\cr}
\eqno(3.6) $$ we get $$ u = {f'(x)\over x\,f'(x)-f(x)}\,. \eqno(3.7)
$$ By a routine computation, $$ {\d u\over\d x} = {-f(x)\,f''(x)
\over(x\,f'(x)-f(x))^2} \asymp 1\qquad \hbox{as } x\to a_i\,, $$ in
view of (3.3). By the second part of (3.6), $$ h(u) = -{u\over
f'(x)} \asymp|x-a_i|^{1-\al_i}\qquad \hbox{as } x\to
a_i\,.\eqno(3.8)  $$ Therefore, using again (3.7) and (3.3), $$
\eq{h'(u) &= {\d\over\d x}\(-{u\over f'(x)}\)\,\({\d u\over\d
x}\)^{-1} \asymp \abs{{\d\over\d x}\((-x\,f'(x)+f(x))^{-1}\)} \cr &
= \abs{x\,f''(x)\over(x\,f'(x)-f(x))^2} \asymp|x-a_i|^{-\al_i}\qquad
\hbox{as } x\to a_i\,.\cr}  $$ Together with (3.8) this implies that
$h(u)\,h'(u)\asymp|x-a_i|^{1-2\al_i}\hbox{ as } x\to a_i$, which
because of $\al_i\le\h$ proves Lemma 3A. \qed

\bsk

{\bf Lemma 3B. } { \it  For a large real parameter $X$ and the
tac-function $H$ defined above, the asymptotics
$$ N(X):=\#\{(m,n)\in\Z\times\Z_{\ge0}:\ H(m,\sqrt{n},0)\le X\ \}
= C\,X^3+\O{X}\eqno(3.9) $$ holds true, with a certain constant
$C>0$. Furthermore, for $0<\delta<1$,
$$ \eq{N_\delta(X):=&\,\#\{(m,n)\in\Z\times\Z_{\ge0}:\ H(m,\sqrt{n},0)\le
X\,,\ \n\le\delta|m|\,\}=\cr =&\,C_\delta\,X^3+\O{X}\,,\cr}
\eqno(3.10)$$ with a positive $C_\delta\ll\delta^2$, the
$O$-constant independent of $\delta$.  As a consequence, for large
$X$ and $0<\omega<1$, $0<\delta<1$, it follows that
$$ \eq{&\#\{(m,n)\in\Z\times\Z_{\ge0}:\ \abs{H(m,\sqrt{n},0)- X} <
\omega\ \} \ll X^2\,\omega+X\,,\cr
&\#\{(m,n)\in\Z\times\Z_{\ge0}:\ \abs{H(m,\sqrt{n},0)- X} <
\omega\,,\ \n\le\delta|m|\,\} \ll X^2\,\omega\,\delta^2+X\,.\cr}
\eqno(3.11) $$} \bsk

{\bf Proof. } Let $\D_+^*$ denote the compact planar domain
bounded by the curve $\c_+^*$ and the $u$-axis. Obviously, $$
\eq{N(X) &= \#\{(m,n)\in\Z\times\Z_{\ge0}:\ (m,\sqrt{n})\in
X\,\D_+^*\ \}= \cr &= \sum_{(1/a_1)X\le m\le
(1/a_2)X}\(1+[X^2\,h^2(m/X)]\)= \cr & = X^2 \sum_{(1/a_1)X\le m\le
(1/a_2)X} h^2(m/X) \ + O(X) = \cr &= X^2
\Int_{(1/a_1)X}^{(1/a_2)X} h^2(u/X)\d u + 2X
\Int_{(1/a_1)X}^{(1/a_2)X} \psi(u)h(u/X)h'(u/X) \d u \ +O(X)\,,
\cr }$$ by the Euler-Mac Laurin formula (see E.~Kr\"atzel [14],
p.~20), with $\psi(u):=u-[u]-\h$. Here the first integral equals
$C\,X$ with $C=\int_{1/a_1}^{1/a_2} h^2(\xi)\d\xi$, which yields
the main term of (3.9). Further, for any interval $[\be_1,
\be_2]\subset\ ]1/a_1, 1/a_2[$, an integration by parts gives
$$ \eq{& 2X \Int_{\be_1X}^{\be_2X} \psi(u)h(u/X)h'(u/X) \d u = \cr =
& 2X\(h(\be_2)h'(\be_2) \psi_1(\be_2 X)-
h(\be_1)h'(\be_1)\psi_1(\be_1X)\)- \cr -& 2 \Int_{\be_1X}^{\be_2X}
\(h(u/X)h''(u/X)+h'^2(u/X)\)\psi_1(u)\d u\,, \cr } \eqno(3.12) $$
where $\psi_1(u):=\int_0^u \psi(v)\d v \ll1$. If $hh''+h'^2=(hh')'$
is bounded on $]1/a_1,1/a_2[$, we simply let $\be_1\to1/a_1$,
$\be_2\to1/a_2$, and obtain the desired bound $O(X)$ for the
remainder. In case that $(hh')'$ is unbounded\fn{3}{By construction,
in particular in view of (3.2), (3.7), (3.8), for $u$ close to
$1/a_i$, $(h(u)h'(u))'$ can be represented as a Laurent series in a
fractional power of $|x-a_i|$. Thus for $u\to1/a_i$,
$|(h(u)h'(u))'|$ either is bounded or tends to $\infty$.} near
$1/a_1$ (say), we choose $\be_1>1/a_1$ such that $(hh')'$ has no
sign change on $]1/a_1,\be_1]$. By the second mean-value theorem and
Lemma 3A,
$$ \Int_{(1/a_1)X}^{\be_1X} \psi(u)h(u/X)h'(u/X) \d u \ll
\sup_{]1/a_1,\be_1]}|hh'|\ll 1\,.  $$ A similar reasoning holds
near $1/a_2$ if necessary. On the remaining interval $[\be_1X,
\be_2X]$, (3.12) readily yields the bound $O(X)$ and thus
completes the proof of (3.9). Quite similarly, $$ \eq{N_\delta(X)
&=  \sum_{(1/a_1)X\le m\le (1/a_2)X}\min\( X^2 h^2(m/X),\delta^2
m^2\) \ + O(X) = \cr &= \Int_{(1/a_1)X}^{(1/a_2)X} \min\( X^2
h^2(u/X),\delta^2 u^2\)\d u \ +\cr &+\ \Int_{(1/a_1)X}^{(1/a_2)X}
\psi(u)\,{\d\over\d u}\(\min\( X^2 h^2(u/X),\delta^2 u^2\)\) \d u
\ +\O{X}\,. \cr }$$ Here the first integral obviously equals
$C_\delta\,X^3$ with $$ C_\delta :=
\Int_{1/a_1}^{1/a_2}\min\(h^2(\xi),
\delta^2\xi^2\)\d\xi\ll\delta^2\,. $$ The remainder integral can
be treated as before, with the bound $O(X)$, since for any
interval $I\subseteq\,]X/a_1,X/a_2[$, $$ \Int_I
\psi(u)\,\delta^2\,u\d u\ll X\,. $$ The deduction of (3.11) from
(3.9), (3.10) is trivial. \qed

\bsk \bsk

{\bf 4.~Asymptotic evaluation of the main terms. } For a large
parameter $t$ it follows, with the definitions of section 3, that
$$ \eq{A_\r(t) &= \sum_{a_1t\le m\le a_2t}\(\sum_{0\le k\le
t^2f^2(m/t)}r(k)\) =\cr &= \pi t^2\sum_{a_1t\le m\le
a_2t}f^2(m/t)+\sum_{a_1t\le m\le a_2t}P(t^2f^2(m/t))\,.\cr}
\eqno(4.1)$$ We proceed to evaluate the first sum on the
right-hand side, postponing the mean-square estimation of the last
one to the next section. By the Euler-Mac Laurin formula,
$$ \eq{\pi t^2\sum_{a_1t\le m\le a_2t}f^2(m/t) &= \pi t^2
\Int_{a_1t}^{a_2t} f^2(\tau/t)\d\tau + \pi t^2
\Int_{a_1t}^{a_2t}\psi(\tau)\,{\d\over\d\tau}(f^2(\tau/t))\d\tau=\cr
&= \vol(\r)t^3 + \pi t^2 \Int_{a_1}^{a_2} \psi(t x)\,{\d\over\d
x}(f^2(x))\d x\,.\cr} \eqno(4.2) $$ By (3.5),
$$ {\d\over\d x}\(f^2(x)\) =
\sum_{j=2}^{2N_1+3} d_{1,j}^* (x-a_1)^{\al_1\,j-1} +
\sum_{j=2}^{2N_2+3} d_{2,j}^* (a_2-x)^{\al_2\,j-1} + \Phi(x)\,,
\eqno(4.3) $$ with $\Phi\in C^1[a_1,a_2]$. Integrating by parts and
using again $\psi_1(u):=\int_0^u \psi(v)\d v \ll1$, we obtain
$$ \eq{t^2 \Int_{a_1}^{a_2} \psi(t x) \Phi(x)\d x =&\
t\,(\psi_1(a_2t)\Phi(a_2)-\psi_1(a_1t)\Phi(a_1))- \cr & - t
\Int_{a_1}^{a_2} \psi_1(t x) \Phi'(x)\d x = O(t)\,.\cr} \eqno(4.4)
$$ The same argument works for $|x-a_i|^{\al_i\,j-1}$
instead of $\Phi(x)$, with $N_i+2\le j\le2N_i+3$, $i\in\{1,2\}$.
In the subsequent analysis we may thus replace the upper summation
limits in the sums from (4.3) by $N_1+1$, resp., $N_2+1$. To deal
with the first one of these remaining sums, we use the Fourier
series
$$ \psi(z)=-{1\over\pi}\sum_{k=1}^\infty{1\over k}\,\sin(2\pi kz)
\qquad(z\notin\Z)$$ and an obvious shift of variable. For
$j=2,\dots,N_1+1$, we conclude that
$$ \Int_{a_1}^{a_2}\psi(tx)(x-a_1)^{\al_1 j-1}\d x =
-{1\over\pi}\sum_{k=1}^\infty {1\over k} \Int_0^{a_2-a_1}x^{\al_1
j-1}\,\sin(2\pi kt(a_1+x))\d x\,. $$ An integration by parts shows
that $$ \Int_{a_2-a_1}^\infty x^{\al_1 j-1}\,\sin(2\pi
kt(a_1+x))\d x =\O{(kt)^{-1}}\,. $$ Further, $$ \Int_0^\infty
x^{\al_1 j-1}\,\sin(2\pi kt(a_1+x))\d x =
\Im\(e(a_1kt)\,(kt)^{-\al_1j}\Int_0^\infty \tau^{\al_1j-1}
e(\tau)\d\tau\) = $$ $$ = \Im\(e(a_1kt)\,(2\pi
kt)^{-\al_1j}\Gamma(\al_1j)\,e(\al_1j/4)\)= \Gamma(\al_1j)\,(2\pi
kt)^{-\al_1j}\,\sin\(2\pi a_1 kt+\h\pi\al_1j\)\,, $$ using
well-known formulas for the last integral (cf., e.g.,
~H.~Rademacher [28], p.~82). Collecting results, we get $$
\eq{&\pi t^2 \Int_{a_1}^{a_2} \psi(t x)\,\(\sum_{j=2}^{N_1+1}
d_{1,j}^* (x-a_1)^{\al_1j-1}\)\d x = \cr & =\ \sum_{j=2}^{N_1+1}
d_{1,j}^*\,t^{2-\al_1j} (2\pi)^{-\al_1j}\,\Gamma(\al_1j)
\sum_{k=1}^\infty k^{-1-\al_1j}\,\sin(-2\pi a_1kt-\h\pi\al_1j) +
\O{t} =\cr & =\ \sum_{j=2}^{N_1+1} d_{1,j}^*\,
\F(-a_1t,\al_1j)\,t^{2-\al_1j} + \O{t}\,, \cr }  $$ where
$\F(\cdot,\cdot)$ has been defined in our Theorem. Quite
similarly, $$ \eq{&\pi t^2 \Int_{a_1}^{a_2} \psi(t
x)\,\(\sum_{j=2}^{N_2+1} d_{2,j}^* (a_2-x)^{\al_2j-1}\)\d x = \cr
& =\ - \sum_{j=2}^{N_2+1} d_{2,j}^*\, \F(a_2
t,\al_2j)\,t^{2-\al_2j} + \O{t}\,. \cr }  $$ Combining the last
two results with (4.2) - (4.4), we finally arrive at $$ \eq{\pi
t^2\sum_{a_1t\le m\le a_2t}f^2(m/t) &= \vol(\r)\,t^3 +
\sum_{j=2}^{N_1+1} d_{1,j}^*\, \F(-a_1t,\al_1j)\,t^{2-\al_1j} -
\cr &- \sum_{j=2}^{N_2+1} d_{2,j}^*\,
\F(a_2t,\al_2j)\,t^{2-\al_2j} + \O{t}\,.\cr } \eqno(4.5) $$

\bsk\bsk

{\bf 5.~Estimating the remainder in mean-square. } It remains to
deal with the last sum in (4.1), i.e., to show that $$ \InT
\(\sum_{a_1 t\le m\le a_2 t}P\(t^2\,f^2(m/t)\)\)^2\d t \ll
T^{3+\ep}\,. $$ For given large $T$ we divide the intervals
$]a_1,\h(a_1+a_2)]$ and \hbox{$]\h(a_1+a_2),a_2]$} into dyadic
subintervals $\J^{(1,r)}=]u^{(1,r+1)},u^{(1,r)}]$,
$\J^{(2,r)}=]u^{(2,r)},u^{(2,r+1)}]$, $0\le r\le R$, where
$u^{(i,r)}:=a_i-(-1)^i\,2^{-r-1}(a_2-a_1)$, and $R$ is chosen such
that the shortest ones of these intervals are of length $\asymp
T^{-1}$. Ignoring the superscripts for short, we write $\J$ for
any of these subintervals, whose number obviously is $O(\log T)$.
Let $\K:=[a_1,a_2]\setminus\bigcup\J$, then $|\K|\asymp T^{-1}$,
and the trivial bound $P(z)\ll\sqrt{z}$ readily implies
$$\InT \(\sum_{m\in t\K}P\(t^2\,f^2(m/t)\)\)^2\d t \ll
T^{3}\,. $$ Thus it suffices to prove that, for each $\J$ and
$\ep>0$,
$$ \InT \(\sum_{m\in t\J}P\(t^2\,f^2(m/t)\)\)^2\d t \ll
T^{3+\ep}\,.\eqno(5.1)
$$ For every $t\in[T,2T]$ and $m\in t\J$, we apply Lemma 1, with
$X=t^2\,f^2(m/t)$ and $Y=T^2$. By (3.3), $X\gg
T^2|\K|^{2\max(\al_i)}\gg T$, hence the condition $\log Y\ll\log
X$ is satisfied. We obtain
$$ \eq{& \sum_{m\in t\J}P\(t^2\,f^2(m/t)\) = \cr & =\
\sum_{m\in t\J}\({\sqrt{t}\over\pi}\,\sqrt{f\({m\over
t}\)}\sum_{1\le n\le T^2}{r(n)\over
n^{3/4}}\,\cos(2\pi\sqrt{n}\,tf(m/t)-3\pi/4) +\O{T^\ep}\) = \cr &
=\ {\sqrt{t}\over\pi} \sum_{1\le n\le T^2}{r(n)\over
n^{3/4}}\left\{ \sum_{m\in t\J}\sqrt{f\({m\over
t}\)}\,\cos(2\pi\sqrt{n}\,tf(m/t)-3\pi/4)\right\} +
\O{T^{1+\ep}}\,.\cr} \eqno(5.2)$$ We shall transform the inner sum
here by means of Lemma 2B, with $$ G(\tau):=\sqrt{f\({\tau/
t}\)}\,,\quad F(\tau):=\sqrt{n}\,tf(\tau/t)\,. $$ To do so we put
$L:=T\,{\rm length}(\J)$, and observe that $|\tau-a t|\asymp L$
for all $\tau\in t\J$, where $a$ is the one of $a_1,a_2$ which is
nearer to $\J$. Hence, in view of (3.3), $$
F''(\tau)=\sqrt{n}\,t^{-1}f''(\tau/t)\asymp
\sqrt{n}\,t^{-1}|\tau/t-a|^{\al-2}\asymp \sqrt{n}\,t^{1-\al}
L^{\al-2} $$ for all $\tau\in t\J$ ($\al$ the appropriate one of
$\al_1,\al_2$), and similarly $$ F^{(j)}(\tau)\ll
\sqrt{n}\,t^{1-\al} L^{\al-j} \for j=3,4\,. $$ Furthermore, by
(3.4), $$
G^{(j)}(\tau))={\d^j\over\d\tau^j}\(\sqrt{f(\tau/t)}\)\ll
t^{-j}\,|\tau/t-a|^{\al/2-j}\ll t^{-\al/2}\,L^{\al/2-j} $$ for
$\tau\in t\J$, $j=0,1,2$. We may thus apply Lemma 2B with the
parameters $$ X:=L\,,\quad
Y:=n^{-1/2}\,t^{\al-1}\,L^{2-\al}\,,\quad Z:=(L/t)^{\al/2}\,.
$$ After a short computation, Lemma 2B yields $$ \eq{&\sum_{m\in t\J}\sqrt{f\({m\over
t}\)}\,\cos(2\pi\sqrt{n}\,tf(m/t)-3\pi/4) = \cr  = {\sqrt{t}\over
n^{1/4}}\, &
\Re\(e(-\h)\sum_{k\in\n\J^*}{\sqrt{f(f^*(k/\sqrt{n}))}\over\abs{f''(f^*(k/\sqrt{n}))}^{1/2}}
\,e\(t(\sqrt{n}f(f^*\kn)-k f^*\kn) \)\) + \cr & +\
\O{n^{-1/4}\,t^{(\al-1)/2}\,L^{1-\al/2}} + \O{\log t}\,,\cr }
\eqno(5.3)
$$ where $f^*$ denotes the inverse function of $f'$ and $\J^*$ the image of the closure
$\bar{\J}$ of $\J$ under $f'$. The contribution of the error terms
here to the whole of (5.2) is $$ \ll  t \sum_{1\le n\le T^2}
{r(n)\over n} + \sqrt{t}\,\log t \sum_{1\le n\le T^2} {r(n)\over
n^{3/4}} \ll T\,\log T\,, $$ hence small enough. We put for short
$\be(k,n):={\sqrt{f(f^*(k/\sqrt{n}))}\over\abs{f''(f^*(k/\sqrt{n}))}^{1/2}}$.
Now $k\in\n\J^*$ implies that $f^*(k/\n)\in\bar{\J}$. Hence, by
(3.3) and the fact that $f''$ is bounded away from zero,
$$ \beta(k,n)\ll {L\over T} \for k\in\n\J^*\,.\eqno(5.4) $$ Furthermore, by the definition of
the tac-function $H$, for all $k\in\Z, n\in\Z^+$, $$
\eq{H(-k,\n,0)&=
\max_{(x,y,0)\in\partial\r}\(-kx+\n\,y\)=\max_{a_1\le x\le
a_2}\(-kx+\n f(x)\)= \cr &= \n f(f^*\kn)-k f^*\kn\,.\cr}  $$ Hence
it will suffice to show that $$ I(\J,T) := \InT t^2 \abs{S(\J,t,T)
}^2 \d t \ll T^{3+\ep}\,, \eqno(5.5)  $$ where $$ S(\J,t,T)
:=\sum_{1\le n\le T^2} {r(n)\over n}\sum_{k\in\n\,\J^*} \beta(k,n)\,
e(t\,H(-k,\n,0))\,.   $$ To simplify the subsequent analysis, we use
a common device involving the Fej\'er kernel $\phi(z):=\sin^2(\pi
z)/(\pi z)^2$. By Jordan's inequality, $\phi(z)\ge4/\pi^2$ for
$\abs{z}\le\h$, and the Fourier transform is simply $$
\fou{\phi}(y)=\Int_{\Ri} \phi(z)\,e(y z) \d y = \max(0,
1-\abs{y})\,. \eqno(5.6)$$

\vbox{Thus $$ I(\J,T)\le4T^3 \Int_{-1/2}^{1/2}
\abs{S(\J,3T/2+Tw,T)}^2 \d w \le $$ $$ \le \pi^2 T^3 \Int_{\Ri}
\phi(w)\,\abs{S(\J,3T/2+Tw,T)}^2 \d w = $$ $$ = \pi^2 T^3
\sum_{1\le n,n'\le T^2} {r(n)r(n')\over
nn'}\sum_{k\in\n\,\J^*\atop k'\in\sqrt{n'}\,\J^*}
\beta(k,n)\beta(k',n')\,\times $$
$$ \times e\({3T\over2}\,\(H(\b m)-H(\b{m'})\)\)\,\fou{\phi}\(T(H(\b m)-H(\b{m'}))\)\,,
$$}

where $\b m:=(-k,\sqrt{n},0)$, $\b{m'}:=(-k',\sqrt{n'},0)$ for
short. Therefore, by (5.4) and (5.6),
$$ I(\J,T)\ll \ T^3\,\lT^2 \sum_{1\le n,n'\le T^2\atop k/\n,\,k'/\ni\in\J^*}
{r(n)r(n')\over nn'}\,\max\(0,1-T|H(\b m)-H(\b{m'})|\)\,. $$ We
observe that $$
\max_{\xi\in\J^*}|\xi|=\max_{x\in\bar{\J}}|f'(x)|\ll\lT^{\al-1}\,,\eqno(5.7)
$$ again by (3.3). We put $\la:=L/T={\rm length}(\J)$ for short.
\def\lT{\la} If $\la$ is sufficiently small (i.e., $\J$ is close to an endpoint
$a_i$), all numbers of $\J^*$ are $\asymp\lT^{\al-1}$. Hence
$k/\n\in\J^*$ implies that $$ \eq{&
k\asymp\n\,\lT^{\al-1}\quad\hbox{ if $\la$ is sufficiently
small,}\cr & k\ll\n\,\lT^{\al-1}\quad\hbox{ always}.\cr} $$
Furthermore, in any case, for $n\ge1$ and $k/\n\in\J^*$,
$$ H(\b m)\asymp|\b m|_\infty\asymp\n\,\lT^{\al-1}\,. $$ Thus $$
\eq{I(\J,T) \ll T^{3+\ep/2}& \lT^{4\al-2} \sum_{1\le n\le T^2\atop
k/\n\,\in\J^*}|\b m|_\infty^{-4}\,\times\cr &\times\,\#\{\b m':\
|H(\b m')-H(\b m)|<T^{-1}\,,\ k'/\ni\,\in\J^*\ \}\,. \cr } $$ By
(3.11) of Lemma 3B,
$$ \#\{\b{m'}=(-k',\ni,0):\ |H(\b{m'})-H(\b m) |<T^{-1},\ \ni\gg|k'|\lT^{1-\al} \,\}\ll
$$ $$
\ll |\b m|_\infty^2\,T^{-1}\,\lT^{2(1-\al)}+|\b m|_\infty\,. $$
Therefore, $$ \eq{I(\J,T) & \ll  T^{3+\ep/2}\sum_{1\le n\le
T^2\atop k\ll\sqrt{n}\,\lT^{\al-1}} \({\lT^{2\al}\over|\b
m|_\infty^{2}\,T}+{\lT^{4\al-2}\over|\b m|_\infty^{3}}\) \ll \cr &
\ll T^{3+\ep/2} \sum_{1\le n\le T^2\atop
k\ll\sqrt{n}\,\lT^{\al-1}} \({\lT^2\over n\,T}+{\lT^{1+\al}\over
n^{3/2}}\)\ll \cr &\ll T^{3+\ep/2} \sum_{1\le n\le T^2}
\({\lT^{1+\al}\over T\,\n} +{\lT^{2\al}\over n}\) \ll T^{3+\ep}\,.
\cr }
$$ This completes the proof of (5.5) and thereby, in view of (4.1)
and (4.5), that of our Theorem. \bsk\msk

%\bye

{\bf6.~Concluding remark. } We indicate briefly how the pointwise
bound (2.3) follows as a by-result from the above analysis.
Estimating the right-hand side of (5.3) trivially, we get, for
$t=T$,
$$ \eq{&\sum_{m\in t\J}\sqrt{f\({m\over
t}\)}\,\cos(2\pi\sqrt{n}\,tf(m/t)-3\pi/4) \ll \cr & \ll
n^{1/4}\,t^{1/2}\,|\J^*| + n^{-1/4}\,t^{(\al-1)/2}\,L^{1-\al/2} +
\log t \ll \cr & \ll n^{1/4}\,t^{3/2-\al}\,L^{\al-1} +
n^{-1/4}\,t^{(\al-1)/2}\,L^{1-\al/2} + \log t  \,,\cr} $$ in view of
(5.7). Using Lemma 1 with $X=t^2\,f^2(m/t)\ll t^2$,
$Y=L^{2-\al}\,t^{\al-1}\ll t$, we obtain as an obvious variant of
(5.2)
$$ \eq{&\sum_{m\in t\J}P\(t^2\,f^2(m/t)\) \ll \cr  \ll
\sqrt{t} & \sum_{1\le n\le Y}  {r(n)\over
n^{3/4}}\,\(n^{1/4}\,t^{3/2-\al}\,L^{\al-1} +
n^{-1/4}\,t^{(\al-1)/2}\,L^{1-\al/2} + \log t\) + {L\,t^{1+\ep}\over
Y^{1/2}} \ll \cr  & \ll t^{2-\al}\,L^{\al-1}\,Y^{1/2} +
t^{\al/2}\,L^{1-\al/2}\,\log t + t^{1/2}\,Y^{1/4}\,\log t +
{L\,t^{1+\ep}\over Y^{1/2}} \ll \cr  & \ll
L^{\al/2}\,t^{3/2-\al/2+\ep}+L^{1-\al/2}\,t^{\al/2} \log
t+L^{1/2-\al/4}\,t^{\al/4+1/4} \log t  \ll t^{3/2+\ep}\,.
 \quad\qed\cr }  $$

\vbox{\vskip 2.5true cm}

{\klein  \parindent=0pt \def\smc{}

\cen{\bf References}  \bsk

[1] {\smc V.~Bentkus \and F.~G\"otze,} On the lattice point
problem for ellipsoids. Acta Arithm. {\bf80}, 101-125 (1997).\ssk

[2] {\smc F.~Chamizo,} Lattice points in bodies of revolution.
Acta Arithm. {\bf85}, 265-277 (1998). \ssk

[3] {\smc F.~G\"otze,} Lattice point problems and values of
quadratic forms. Invent.~Math. {\bf157}, 195--226 (2004).\ssk

[4] {\smc K.~Haberland}, \"{U}ber die Anzahl der Gitterpunkte in
konvexen Gebieten. Preprint FSU Jena 1993 (unpublished).\ssk

[5] {\smc E.~Hlawka,} \"Uber Integrale auf konvexen K\"orpern I.
Monatsh.~f.~Math. {\bf 54}, 1-36 (1950). \ssk

[6] {\smc E.~Hlawka}, \"Uber Integrale auf konvexen K\"orpern II.
Monatsh.~f.~Math. {\bf 54}, 81--99 (1950). \ssk

[7] {\smc M.N.~Huxley}, {Area, lattice points, and exponential
sums.} LMS Monographs, New Ser. {\bf 13}, Oxford 1996.  \ssk

[8] {\smc A.~Iosevich, E.~Sawyer, \and A.~Seeger,} Mean square
discrepancy bounds for the number of lattice points in large
convex bodies. J.~Anal.~Math. {\bf87}, 209-230 (2002). \ssk

[9] {\smc A.~Ivi\'c}, {The Laplace transform of the square in the
circle and divisor problems,} Stud.~Sci.~Math. Hung. {\bf32},
181-205 (1996). \ssk

[10] {\smc A.~Ivi\'c, E.~Kr\"{a}tzel, M.~K\"{u}hleitner, \and W.G.~Nowak,}
Lattice points in large regions and related arithmetic functions:
Recent developments in a very classic topic. Proceedings Conf.~on
Elementary and Analytic Number Theory ELAZ'04, held in Mainz, May
24-28, W.~Schwarz and J.~Steuding eds., Franz Steiner Verlag 2006,
pp. 89-128.\ssk

[11] {\smc V.~Jarnik}, \"Uber die Mittelwerts\"atze der
Gitterpunktlehre, V.~Abh. Cas.~mat.~fys. {\bf69}, 148--174
(1940).\ssk

[12] {\smc A.A.~Karatsuba \and S.M.~Voronin}, The Riemann
zeta-function. Berlin 1992. \ssk

[13] {\smc D.G.~Kendall}, On the number of lattice points inside a
random oval. Quart.~J.~Math.~(Oxford) {\bf 19}, 1-26 (1948).  \ssk

[14] {\smc E.~Kr\"atzel,} Lattice points. Berlin 1988. \ssk

[15] {\smc E.~Kr\"atzel,} Analytische Funktionen in der
Zahlentheorie. Stuttgart-Leipzig-Wiesbaden 2000.  \ssk

[16] {\smc E.~Kr\"atzel,} Lattice points in three-dimensional large
convex bodies. Math.~Nachr. {\bf212}, 77--90 (2000). \ssk

[17] {\smc E.~Kr\"atzel,} Lattice points in three-dimensional convex
bodies with points of Gaussian curvature zero at the boundary.
Monatsh.~Math. {\bf137}, 197--211 (2002).  \ssk

[18] {\smc E.~Kr\"atzel,} Lattice points in some special
three-dimensional convex bodies with points of Gaussian curvature
zero at the boundary. Comment.~Math.~Univ.~Carolinae {\bf43},
755-771 (2002). \ssk

[19] {\smc E.~Kr\"atzel \and W.G.~Nowak}, Eine explizite
Absch\"{a}tzung f\"{u}r die Gitter-Diskrepanz von Rotations\-ellipsoiden.
Monatsh.~Math., to appear.  \ssk

[20] {\smc M.~K\"uhleitner,} Lattice points in bodies of
revolution in $\Ri^3$: an $\Omega_-$-estimate for the error term.
Arch.~Math.~(Basel) {\bf 74}, 234-240 (2000). \ssk

[21] {\smc M.~K\"uhleitner \and W.G.~Nowak,} The lattice point
discrepancy of a body of revolution: Improving the lower bound by
Soundararajan's method. Arch.~Math. (Basel) {\bf83}, 208--216
(2004). \ssk

[22] {\smc Y.-K.~Lau,} On the mean square formula of the error term
for a class of arithmetical functions. Monatsh.~Math. {\bf128},
111-129 (1999). \ssk

[23] {\smc W.~M\"uller}, On the average order of the lattice rest of
a convex body. Acta Arithm. {\bf 80}, 89--100 (1997).  \ssk

[24] {\smc W.~M\"uller}, Lattice points in large convex bodies.
Monatsh.~Math. {\bf128}, 315--330 (1999).  \ssk

[25] {\smc W.G.~Nowak}, On the average order of the lattice rest of
a convex planar domain. Proc.~Cambridge Phil.~Soc. {\bf98}, 1--4
(1985).\ssk

[26] {\smc W.G.~Nowak}, On the mean lattice point discrepancy of a
convex disc. Arch.~Math.~(Basel) {\bf78}, 241--248 (2002).\ssk

[27] {\smc M.~Peter}, Lattice points in convex bodies with planar
points on the boundary. Monatsh.~Math. {\bf135}, 37--57 (2002). \ssk

[28] {\smc H.~Rademacher}, Topics in analytic number theory.
Berlin-Heidelberg-New York 1973.

\vbox{\vskip 1true cm}

\parindent=1.5true cm

\vbox{Werner Georg Nowak \ssk

Institute of Mathematics

Department of Integrative Biology

Universit\"at f\"ur Bodenkultur Wien

Gregor Mendel-Stra\ss e 33

A-1180 Wien, \"Osterreich \ssk

E-mail: {\tt nowak@boku.ac.at} \ssk

}}

\bye